\documentclass{article}

\usepackage{amsmath,amsfonts,amssymb,amsthm}

\begin{document}

\begin{center}
AN EFFICIENT METHODOLOGY FOR THE ANALYSIS AND MODELING OF COMPUTER EXPERIMENTS WITH LARGE NUMBER OF INPUTS

Bertrand Iooss$^{12}$ and Amandine Marrel$^3$


$^1$EDF R\&D \\
  6 Quai Watier, 78401 Chatou, France\\
  e-mail: bertrand.iooss@edf.fr \\
  $^2$
  Institut de Math\'ematiques de Toulouse\\
  31062 Toulouse, France \\
  $^3$
  CEA, DEN, DER \\
  13108 Saint-Paul-lez-Durance, France\\
  e-mail: amandine.marrel@cea.fr

\end{center}

\begin{abstract}

Complex computer codes are often too time expensive to be directly used to perform uncertainty, sensitivity, optimization and robustness analyses. A widely accepted method to circumvent this problem consists in replacing cpu-time expensive computer models by cpu inexpensive mathematical functions, called metamodels. For example, the Gaussian process (Gp) model has shown strong capabilities to solve practical problems, often involving several interlinked issues. However, in case of high dimensional experiments (with typically several tens of inputs), the Gp metamodel building process remains difficult, even unfeasible, and application of variable selection techniques cannot be avoided.  
In this paper, we present a general methodology allowing to build a Gp metamodel with large number of inputs in a very efficient manner. While our work focused on the Gp metamodel, its principles are fully generic and can be applied to any types of metamodel. The objective is twofold: estimating from a minimal number of computer experiments a highly predictive metamodel. 
This methodology is successfully applied on an industrial computer code.
\end{abstract}

Keywords: Computer experiments, Uncertainty Quantification,Sensitivity Analysis, Metamodel, Gaussian process

\section{INTRODUCTION}

Quantitative assessment of the uncertainties tainting the results of computer simulations is nowadays a major topic of interest in both industrial and scientific communities.
One of the key issues in such studies is to get information about the output when the numerical simulations are expensive to run. 
For example, in nuclear engineering problems, one often faces up with cpu time consuming numerical models and, in such cases, uncertainty propagation, sensitivity analysis, optimization processing and system robustness analysis become difficult tasks using such models.
In order to circumvent this problem, a widely accepted method consists in replacing cpu-time expensive computer models by cpu inexpensive mathematical functions, called metamodels \cite{fanli06}.
This solution has been applied extensively and has shown its relevance especially when simulated phenomena are related to a small number of random input variables (see \cite{forsob08} for example).

However, in case of high dimensional numerical experiments (with typically several tens of inputs), depending on the complexity of the underlying numerical model, the metamodel building process remains difficult, even unfeasible.
For example, the Gaussian process (Gp) model \cite{sanwil03} which has shown strong capabilities to solve practical problems, has some caveats when dealing with high dimensional problems. 
The main difficulty relies on the estimation of Gp hyperparameters.
Manipulating pre-defined or well-adapted Gp kernels (as in  \cite{muerou12,durgin13}) is a current research way, while coupling the estimation procedure with variable selection techniques has been proposed by several authors \cite{welbuc92,marioo08,woolew17}.

In this paper, we pursue the effort on the latter technique by proposing a more rigorous and robust method for building a Gp metamodel with a high-dimensional vector of inputs.
First, we clarify the sequence of the different steps of the methodology, while updating their technical core with more relevant statistical techniques.
For example, the screening step is raised by the use of recent and powerful techniques in terms of variable selection using a small number of model runs.
Second, contrary to the previous works, we do not remove the non-selected inputs from the Gp model, keeping the uncertainty caused by the dimension reduction by using the joint metamodel technique \cite{marioo12}.
The integration of this residual uncertainty is important in terms of robustness of subsequent safety studies.

The next section of this paper presents our general methodology.
The third, fourth and fifth sections are devoted to a detailed explanation of each of its steps.
The last section shows an application of this work on a thermal-hydraulic calculation case simulating accidental scenario in a nuclear reactor.
It also gives some prospects of this work.

\section{GENERAL METHODOLOGY}

The system under study is denoted
\begin{equation}
Y=g\left(X_1,\ldots,X_d\right)
\end{equation}
where $g(\cdot)$ is the numerical model (also called the computer code), whose output $Y$ and input parameters $X_1,\ldots,X_d$ belong to some measurable spaces $\mathcal{Y}$ and $\mathcal{X}_1, \ldots, \mathcal{X}_d$ respectively. 
$\mathbf{X}=\left(X_1,\ldots,X_d\right)$ is the input vector and we suppose that $\mathcal{X}=\prod_{k=1}^d\mathcal{X}_k \subset  \mathbb{R}^d $ and $\mathcal{Y}\subset \mathbb{R}$. 
For a given value of the vector of inputs $\mathbf{x} = \left(x_1,\ldots,x_d\right) \in \mathbb{R}^d$, a simulation run of the code yields an observed value $y = g(\mathbf{x})$. 

Our approach consists in four steps:
\begin{enumerate}

\item \textbf{Step 1: Initial experimental design.} Once the uncertain input variables of the numerical model $g$ and their variation domain identified, a design of $n$ experiments is firstly performed and yields $n$ model output values. To constitute this learning sample, we use a space-filling design (SFD) of experiments, providing a full coverage of the high-dimensional input space.

\item \textbf{Step 2: Screening.} From the learning sample, a screening technique is performed in order to identify the primary influential inputs (PII) on the model output variability. It has been recently shown that screening based on dependence measures \cite{dav15,delmar16b} or on derivative-based global sensitivity measures \cite{kucioo17,roubar16} are very efficient methods which can be directly applied on a SFD. 
One of their great interest is that, additionally to their screening job, the sensitivity indices they provide can be quantitatively interpreted.
From these screening results, the inputs are then ordered by decreasing PII, for the purpose of the metamodeling step.

\item \textbf{Step 3: Joint metamodeling and metamodel validation.} The sorted inputs are successively included in the group of explanatory inputs while the other inputs are considered as a global stochastic (\textit{i.e. unknown}) input and a joint Gp metamodel is built. At each iteration, a first Gp model, only depending on the explanatory inputs, is built to approximate the mean component of the metamodel. The residual effect of the other inputs is captured using a second Gp model which approximates the variance component as a function of the explanatory inputs. For this, a joint metamodeling procedure is used, as proposed by \cite{marioo12}. Moreover, in order to deal with the large number of inputs, the optimization process, which is required to estimate the hyperparameters of the Gp covariances, uses as a starting point the values estimated at the previous iteration.

The accuracy and prediction capabilities of the metamodel are controlled on a test sample or by cross-validation. 
\end{enumerate}

All these steps are described in the next subsections.
The obtained metamodel, which requires a negligible calculation time, can then be used to perform global sensitivity analysis, uncertainty propagation (for example through Monte-Carlo simulations) or optimization processing.

\section{STEP 1: INITIAL DESIGN OF EXPERIMENTS}\label{sec:SFD}

The objective of the initial sampling step is to investigate the whole variation domain of the uncertain parameters in order to fit a predictive metamodel which approximates as accurately as possible the code in the whole domain of variation of the uncertain parameter, independently from their probabilistic distributions. 
For this, we use a space-filling design (SFD) of a certain number $n$ of experiments, providing a full coverage of the high-dimensional input space \cite{fanli06}.
This design enables to investigate the domain of variation of the uncertain parameters and provides a learning sample.

Mathematically, this corresponds to the sample $\left\{ \mathbf{x}^{(1)},\ldots,\mathbf{x}^{(n)}\right\}$ which is performed on the model $g$.
This yields $n$ model output values denoted $\left\{ y^{(1)},\ldots,y^{(n)} \right\}$ with $y^{(i)} = g(\mathbf{x}^{(i)})$. 
The obtained learning sample is denoted $\left( X_s, Y_s\right)$ with $X_s = \left[{\mathbf{x}^{(1)}}^T,\ldots,{\mathbf{x}^{(n)}}^T\right]^T$ and $Y_s = \left[y^{(1)},\ldots,y^{(n)}\right]^T$. 
The goal is to build an approximating model of $g$ using the $n$-sample $\left( X_s, Y_s\right)$. 

The number $n$ of simulations is a compromise between the CPU time required for each simulation and the number of input parameters.
Some thumb rules propose to choose $n$ at least as large as $10$ times the dimension $d$ of the input vector \cite{loesac09,marioo08}.

For the SFD type, a Latin Hypercube Sample (LHS) with optimal space-filling and good projection properties \cite{woolew17} would be well adapted.
In particular, \cite{fanli06,damcou13} have shown the importance of ensuring good low-order sub-projection properties.
Maximum projection designs \cite{josgul15} or low-centered $L^2$ discrepancy LHS \cite{jinche05} are then particularly well-suited.

{\it Remark:
Note that the input values are sampled uniformly, considering only their variation ranges and not their initial probability distributions. 
Indeed, our aim is to build a metamodel for a multi-objective purpose (sensitivity analysis, uncertainty propagation, etc...).
The input probability distributions will then be used in the sensitivity analysis or uncertainty propagation studies. 
}

\section{STEP 2: INITIAL SCREENING}

From the learning sample, an initial screening is performed in order to identify the PII and sort them by decreasing order of influence. 
For this, two possibilities are proposed: one based on dependence measures and another based on derivative-based global sensitivity measures.

\subsection{Screening based on dependence measure}

\cite{dav15} and more recently \cite{delmar16b} have proposed to use dependence measures for screening purpose, by applying them directly on a SFD. 
These sensitivity indices are not the classical ones variance-based measures (see \cite{ioolem15} for a global review).
They consider higher order information about the output behavior in order to provide more detailed information. 
Among them, the Hilbert-Schmidt independence criterion (HSIC) introduced by \cite{gretbou05} builds upon kernel-based approaches
for detecting dependence, and more particularly on cross-covariance operators in reproducing kernel Hilbert spaces (RKHS).

If we consider two RKHS $\mathcal{F}_k$ and $\mathcal{G}$ of functions $\mathcal{X}_k \rightarrow \mathbb{R}$ and $\mathcal{Y} \rightarrow \mathbb{R}$ respectively, the crossed-covariance $C_{X_k, Y}$ operator associated to the joint distribution of $\left(X_k,Y\right)$ is the linear operator defined for every $f_{X_k} \in \mathcal{F}_k$ and $g_Y \in \mathcal{G}$ by: 
\begin{equation}
\langle f_{X_k},C_{X_k, Y} g_{Y} \rangle_{\mathcal{F}_k}=\text{Cov}\left(f_{X_k},g_{Y}\right).
\end{equation}
$C_{X_k, Y}$ generalizes the covariance matrix by representing higher order correlations between $X_k$ and $Y$ through nonlinear kernels. The HSIC criterion is then defined by the Hilbert-Schmidt norm of the
cross-covariance operator:
\begin{equation}
\text{HSIC}(X_k,Y)_{\mathcal{F}_k,\mathcal{G}} = \|C_k\|_{HS}^2.
\end{equation}
From this, \cite{dav15} introduces a normalized version of the HSIC which provides a sensitivity index of $X_k$:
\begin{equation}\label{eq:HSICRk}
 R^2_{\text{HSIC},k}=\frac{\text{HSIC}(X_k,Y)}{\sqrt{\text{HSIC}(X_k,X_k)\text{HSIC}(Y,Y)}}.
\end{equation}
\cite{gretbou05} also propose a Monte Carlo estimator of $\text{HSIC}(X_k,Y)$ and a plug-in estimator can be deduced for $R^2_{\text{HSIC},k}$. Note that Gaussian kernel functions with empirical estimations of the variance parameter are used in our application  (see \cite{gretbou05} for details).

Then, from the estimated $R^2_{\text{HSIC}}$, independence tests are performed for a screening purpose. The objective is to separate the inputs into two sub-groups, the significant ones and the non-significant ones.
For a given input $X_k$, it aims at testing the null hypothesis ``$\mathcal{H}_0^{(k)}$: $X_k$ and $Y$ are independent'', against its alternative ``$\mathcal{H}_1^{(k)}$: $X_k$ and $Y$ are dependent''. The significance level\footnote{The significance level of a statistical hypothesis test is the rate of the type I error which corresponds to the rejection of the null hypothesis $\mathcal{H}_0$ when it is true.} of these tests is hereinafter noted $\alpha$. Several statistical hypothesis tests are available: asymptotic versions, spectral extensions and bootstrap versions for non-asymptotic case. All these tests are described and compared in \cite{delmar16b}; a guidance to use them for a screening purpose is also proposed.
At the end of the screening step, the inputs selected as significant are also ordered by decreasing $R^2_{\text{HSIC}}$. This order will be used for the sequential metamodel building in step 3.

\subsection{Screening based on derivative-based global sensitivity measure}\label{sec:dgsm}

The so-called Derivative-based Global Sensitivity Measures (DGSM) consist in integrating the square derivatives of the model output (with respect to each of the model input) over the domain of the inputs.
This kind of indices have been shown to be easily and efficiently estimated by sampling techniques (as Monte Carlo or quasi-Monte Carlo).
Several authors have shown the interest of DGSM as a screening technique (see \cite{kucioo17} for a review).
Indeed, the DGSM interpretation is made easy due to its inequality links with variance-based sensitivity indices, which are easily interpretable \cite{ioolem15}.
Multiplied by an optimal Poincar\'e constant, DGSM is a narrow upper bound of the total Sobol' index \cite{roubar16}, whatever the input probability distribution.

One of the main issue for this technique in practical situations is to efficiently estimate the model derivatives as the standard practice based on finite-differences is relatively costly.
Indeed, its cost linearly depends on the number of inputs as most of the sensitivity analysis techniques \cite{ioolem15}.
However, if the reverse (adjoint) mode of the numerical model is available, computing all partial derivatives of the model output has a cost independent on the number of input variables.
In this case, the screening step can be performed with a reasonable cpu time cost (with a sample of $100$ runs of the adjoint model typically) and is therefore possible even for large-dimensional model.
This potentiality has been recently applied in \cite{petgoe16} which studies a model with $40$ inputs and uses automatic differentiation in order to obtain the adjoint model (which has a cost of two times the direct model). 
On this example, \cite{roubar16} have shown the relevance of DGSM for a screening purpose, which also provides a quantitative interpretation of sensitivity indices.

\section{STEP 3: JOINT GP METAMODEL WITH SEQUENTIAL BUILDING PROCESS}\label{seq:JM_building}

Among all the metamodel-based solutions (polynomials, splines, neural networks, etc.), we focus our attention on the Gaussian process (Gp) regression, which extends the kriging principles of geostatistics to computer experiments by considering the correlation between two responses of a computer code depending on the distance between input variables. 
The Gp-based metamodel presents some real advantages compared to other metamodels: exact interpolation property, simple analytical formulations of the predictor, availability of the mean squared error of the predictions and the proved efficiency of the model \cite{sanwil03}.
 
However, for its application to complex industrial problems, developing a robust implementation methodology is required. Indeed, fitting a Gp model implies the estimation of several hyperparameters involved in the covariance function. In complex situations (e.g. large number of inputs), some difficulties can arise from the parameter estimation procedure (instability, high number of hyperparameters, see \cite{marioo08} for example). To tackle this issue, we propose a progressive estimation procedure which combines the result of the previous screening step and a joint Gp approach \cite{marioo12}.

\subsection{Successive inclusion of explanatory variables}\label{methodo:sucess}

At the end of the screening step, the inputs selected as significant are ordered by decreasing influence. The sorted inputs thus obtained are successively included in the group of explanatory inputs. At the $j^{th}$ iteration, only the $j$ first sorted inputs are considered as explanatory input variables while all the remaining inputs are included in a single macro-parameter.
This macro-parameter is considered as an uncontrollable parameter (i.e. a stochastic parameter, notion detailed in section \ref{jointGp}). 

At this stage, a joint Gp metamodel is then built with the $j$ explanatory inputs, following the procedure described in \cite{marioo12} and summarized in the next subsection. A numerical optimization is performed to estimate the parameters of the joint metamodel (covariance and variance parameters). In order to improve the robustness of the optimization process, the estimated hyperparameters obtained at the $(j-1)^{th}$ iteration are used, as starting points for the optimization algorithm. This procedure is repeated until the inclusion of all the significant input variables.
Note that this sequential process is directly adapted from the one proposed by \cite{marioo08}.

\subsection{Joint Gp metamodel}\label{jointGp}

In the framework of stochastic computer codes, \cite{zabdej98} proposed to model the mean and dispersion of the code output by two interlinked Generalized Linear Models (GLM), called ``joint GLM''.
\cite{marioo12} extends this approach to several nonparametric models and obtains the best results with two interlinked Gp models, called ``joint Gp''. 
In this case, the stochastic input is considered as an uncontrollable parameter denoted $\mathbf{X}_\varepsilon$ (i.e. governed by a seed variable). 

We extend this approach to a group of non-explanatory variables. More precisely, the input variables $\mathbf{X}=(X_1,\ldots,X_d)$ are divided in two subgroups: the explanatory ones denoted $\mathbf{X_{exp}}$ and the others denoted $\mathbf{X}_\varepsilon$.
The output is thus defined by $ y = g(\mathbf{X_{exp}},\mathbf{X}_\varepsilon)$.
Under this hypothesis, the joint metamodeling approach yields building two metamodels, one for the mean $Y_m$ and another for the dispersion component $Y_d$:
\begin{equation}\label{eqYm}
  Y_m(\mathbf{X_{exp}}) = \mathbb{E}(Y|\mathbf{X_{exp}}) 
\end{equation}
\begin{equation}\label{eqYd}
  Y_d(\mathbf{X_{exp}}) = \mbox{Var}(Y|\mathbf{X_{exp}}) = \mathbb{E}\left[ (Y-Y_m(\mathbf{X_{exp}}))^2 |\mathbf{X_{exp}} \right]. 
\end{equation}

To fit these mean and dispersion components, we propose to use the methodology proposed by \cite{marioo12}. First, an initial Gp denoted $Gp_{m,1}$ is estimated for the mean component with homoscedastic nugget effect. A nugget effect is required to relax the interpolation property of the Gp metamodel, which would yield zero residuals for the whole learning sample. Then, a second Gp, denoted $Gp_{v,1}$, is built for the dispersion component with, here also, an homoscedastic nugget effect. $Gp_{v,1}$ is fitted on the squared residuals from the predictor of $Gp_{m,1}$. Its predictor is considered as an estimator of the dispersion component. The predictor of $Gp_{v,1}$ provides an estimation of the dispersion at each point. It is thus considered as the value of the heteroscedastic nugget effect: the homoscedastic hypothesis is removed. A new Gp, $Gp_{m,2}$, is fitted on data, with the estimated heteroscedastic nugget. Finally, the Gp on the dispersion component is updated from $Gp_{m,2}$ following the same methodology as the one $Gp_{v,1}$.

{\it Remark:
Note that some parametric choices are made for all the Gp metamodels: a constant trend and a Mat\'{e}rn stationary anisotropic covariance are chosen. All the hyperparameters (covariance parameters) and the nugget effect (when homoscedastic hypothesis is done) are estimated by maximum likelihood optimization process.
}

\subsection{Assessment of metamodel accuracy}

To evaluate the accuracy of the metamodel, we use the predictivity coefficient $Q^2$: 
\begin{equation}
Q^2=1-\frac{\sum_{i=1}^{n_{\text{test}}}\left(y^{(i)}-\hat{y}^{(i)}\right)^2}{\sum_{i=1}^{n_{\text{test}}} \left( y^{(i)} - \frac{1}{n_{\text{test}}}  \sum_{i=1}^{n_{\text{test}}}  y^{(i)} \right)^2}
\end{equation}
where $(x^{(i)})_{1\leq i \leq n_{\text{test}}}$ is a test sample, $(y^{(i)})_{1\leq i \leq n_{\text{test}}}$  are the corresponding observed outputs and $(\hat{y}^{(i)})_{1\leq i \leq n_{\text{test}}}$ are the metamodel predictions. $Q^2$ corresponds to the coefficient of determination in prediction and can be computed on a test sample independent from the learning sample or by cross-validation on the learning sample. The closer to one the $Q^2$, the better the accuracy of the metamodel. Note that, in our sequential building process (cf. Section \ref{methodo:sucess}), the $Q^2$ coefficient metamodel is computed at each iteration.

In the case where the model provides the adjoint code (see Section \ref{sec:dgsm}), the gradient evaluations could be integrated in the metamodel building. For this, the co-kriging principle could be adapted to the joint metamodel approach.

\section{APPLICATION TO A THERMAL-HYDRAULIC COMPUTER CODE}

\subsection{Description of the use-case}

Our use-case consists in thermal-hydraulic computer experiments, typically used in support of regulatory work and nuclear power plant design and operation. 
Indeed, some safety analysis considers the so-called ``Loss Of Coolant Accident'' (LOCA), which takes into account a double-ended guillotine break with a specific size piping rupture.
It is modeled with code CATHARE 2.V2.5 which simulated the thermalhydraulic responses during a LOCA in a Pressurized water Reactor \cite{mazvac16}. 

In this use-case, $27$ scalar input variables of CATHARE are uncertain.
In our problem, they are defined by their minimum and maximum.
They correspond to various system parameters as initial conditions, boundary conditions, some critical flowrates, interfacial friction coefficients, condensation coefficients, \ldots
The output variable of interest is a single scalar which is the maximal peak cladding temperature during the accident transient.
Our objective with this use-case is to provide a good metamodel to the safety engineers.
Indeed, the cpu-time cost of this computer code is too important to develop all the statistical analysis required in a safety study only using direct calculations of the computer code.
A metamodel would allow to develop more complete and robust demonstration.

$1000$ CATHARE simulations of this test case have been provided following a space-filling LHS with good projection properties (see Section \ref{sec:SFD}) as the design of experiments.
In this test case, the adjoint model is not available and the derivatives of the model output are therefore not computed because of their costs.
The screening step will then be based on HSIC, obtained from the inputs-output sample.

In order to test it, our overall methodology is applied with different sizes of the learning sample: $n=200$, $400$, $600$ and $800$. In each case, the remaining simulations (from the $1000$ we have) are used as a test sample in order to compute the metamodel predictivity. 

\subsection{Screening step with HSIC}

The normalized HSIC coefficients are computed for the different learning sample sizes. Similar results are obtained. Four variables are identified as the most influential: $X_{10}$ (HSIC $\approx 30\%$), $X_{12}$ and $X_{13}$ (HSIC $\approx 14\%$) and $X_{22}$ (HSIC $\approx 9\%$). $X_{14}$, $X_{15}$ and $X_{2}$ have also a significant but lower influence (HSIC around 5\%). Thus, statistical significance tests (asymptotic version with $\alpha = 10 \%$) have selected these 7 inputs. The estimated HSIC and the results of significant tests are relatively stable and independent from the learning sample size, only one or two additional variables with a very low HSIC ($< 2 \%$) are selected for the smallest sample size. 
This confirms the robustness of the HSIC indices and the associated significance tests for qualitative sorting and screening purpose. 

For each learning sample size, the significant inputs are considered as the explanatory variables in the joint metamodel and will be successively included in the building process. The other variables are joined in the uncontrollable parameter. 

\subsection{Joint Gp}\label{appli:GP}

From the HSIC-based screening results, the joint Gp metamodel is built following the sequential process described in Section \ref{seq:JM_building}. 
The simple Gp metamodel with all the $27$ inputs as explanatory variables is also built, without any sequential approach. To assess the accuracy of the different metamodels, the predictivity coefficient $Q^2$ is computed by cross-validation (leave-one-out process) and on the test sample composed of the remaning simulations. The $Q^2$ obtained with the different metamodels are presented in Table \ref{tab:predictivity}. Note that the same optimizer is used to estimate the hyperparameters by maximum likelihood, in order to allow for a fair comparison.

\begin{table}[!ht]
	\centering
	\begin{tabular}{|c||c|c||c|c|}
		\hline
		& \multicolumn{2}{c||}{\textbf{Joint Gp with sequential approach}} & \multicolumn{2}{c|}{\textbf{Simple Gp without sequential approach}}\\
		 & $Q^2$  & $Q^2$  & $Q^2$ & $Q^2$\\
		 & on test sample & by cross-validation & on test sample & by cross-validation\\
		\hline
		$n=200$ & 0.82 & 0.81 & 0.75 & 0.78 \\
		$n=400$ & 0.82 & 0.85 & 0.78 & 0.85 \\
		$n=600$ & 0.86 & 0.89 & 0.83 & 0.86 \\
		$n=600$ & 0.87 & 0.88 & 0.82 & 0.84 \\
	\hline
	\end{tabular}
	\caption{Comparison of Gp metamodel predictivity for different sizes $n$ of learning sample and different building processes.}
	\label{tab:predictivity}
\end{table}	

The joint Gp with a sequential building process outperforms the simple Gp directly built with the 27 input variables, especially for the lower learning sample sizes. On average, the $Q^2$ is improved between 3\% to 9\%. Thus, the proposed methodology allows a more robust metamodel building with a high-dimensional vector of inputs, even with small sample sizes. Moreover, even if it is not used and illustrated in this application, the dispersion component of the joint metamodel takes into account the uncertainty due to the non-significant inputs. This residual uncertainty, although low, is kept by using the joint metamodel technique.
It appears in the mean squared error of the metamodel predictions and could be integrated in subsequent sensitivity or uncertainty propagation studies.
	
\subsection{Work continuation and prospects}

Using the fitted joint Gp metamodel, several statistical analysis, not feasible with the numerical model due to its computational cost, are now accessible.
First, variance-based sensitivity analysis using Sobol' indices can be fully expressed using a Gp model \cite{marioo08,legcan14}.
This would provide a fine determination of the critical parameters whose uncertainty has to be reduced.

Second, we are particularly interested by the estimation of high quantile (at the order of $95\%$ to $99\%$) of the model output temperature.
In nuclear safety, methods of conservative computation of quantiles \cite{nutwal04} have been largely studied.
However, several complementary information are often useful and are not accessible in a high-dimensional context.
Then, we expect the Gp metamodel can help to access this information.
For instance, quantile-based sensitivity analysis \cite{brofor17} and quantile robustness analysis (using the sensitivity indices called PLI \cite{suebou16}) are fully devoted to quantile.
Their relevance to support safety analysis seems promising.


%


\section{Acknowledgments}

We are grateful to Henri Geiser and Thibault Delage who performed the computations of the CATHARE code.

\bibliographystyle{plain}

\end{document}